\documentclass[12pt]{article}
\usepackage{amssymb, url}
\usepackage{amsthm}
\title{On Permutations Avoiding Arithmetic Progressions}
\author{Timothy D. LeSaulnier $\,$ and $\,$ Sujith Vijay \\\\
{\em{\normalsize{Department of Mathematics}}} \\ 
{\em{\normalsize{University of Illinois at Urbana-Champaign}}} \\ 
{\em{\normalsize{Urbana, IL 61801, USA.}}}}

\newtheorem{theorem}{Theorem}

\newcommand{\eop}{\nobreak \ifvmode \relax \else
\ifdim\lastskip<1.5em \hskip-\lastskip
\hskip1.5em plus0em minus0.5em \fi \nobreak
\vrule height0.75em width0.5em depth0.25em\fi}

\begin{document}

\date{}
\maketitle

Let $M(n)$ denote the number of permutations of $\{1,2,\ldots,n\}$ that do 
not contain a $3$-term arithmetic progression as a subsequence. For 
example, $M(4)=10$, corresponding to the permutations $(1243), \, 
(1324), \, (2143), \, (2413), \, (4213)$ and their reversals. In 1977, 
Davis, Entringer, Graham and Simmons \cite{Dav77} established the 
following bounds on $M(n)$: $$2^{n-1} \le M(n) \le \lfloor (n+1)/2 \rfloor 
\, ! \lceil (n+1)/2 \rceil \, !$$ These bounds were recently improved by 
Sharma \cite{Sha09}, who showed that $$M(n) \le \frac{2.7^n}{21} \mbox{ 
for } n \ge 11$$ and that $$\lim_{n \to \infty} \frac{M(n)}{2^n \: 
n^k} = \infty \mbox { for any fixed } k.$$

In \cite{Sha09} the question whether $\lim_{_{{}_{\! \! \! \! \! \! \! \! 
\! \! \! \! \! \!  n \to \infty}}} \frac{M(n+1)}{M(n)} = 2$ was attributed 
to the Editor of the Problem Section of the American Mathematical Monthly 
(where the function $M(n)$ made its earliest known appearance, in 1975), 
and was mentioned as still open. We begin with an observation that settles 
this question in the negative. Indeed, we establish the following stronger 
lower bound on $M(n)$. \\

\theorem $M(n) \ge (1/2)c^n$ for $n \ge 8$, where $c=(2132)^{1/10}=2.152...$. 

\proof The following inequalities were proved in \cite{Dav77} to show that 
$M(n) \ge 2^{n-1}$: $$M(2n) \ge 2[M(n)]^2 \; ; \; M(2n+1) \ge 
2M(n)M(n+1)$$ These recurrences follow from the observation that if 
$\sigma_1$ and $\sigma_2$ are $3AP$-free permutations of 
$\{2,4,\ldots,2n\}$ and $\{1,3,\ldots,2n-1\}$, concatenating them in 
either order yields $3AP$-free permutations $\sigma_1 \sigma_2$ and 
$\sigma_2 \sigma_1$ of $\{1,2,\ldots,2n\}$, since the first and third 
terms of any 3AP have the same parity. Note that these recurrences imply 
the stronger lower bound $M(n) \ge (1/2)c^n$ for $n \ge 8$, where 
$c=(2M(10))^{1/10}=2.152...$. Since $M(8)=282, \, M(9)=496, \, M(10)=1066, 
\, M(11)=2460, \, M(12)=6128, \, M(13)=12840, \, M(14)=29380$ and 
$M(15)=73904$ (see \cite{Dav77}), the inequality holds for $8 \le n \le 
15$. We can now use induction on $k$ to show that it also holds for 
$2^k \le n < 2^{k+1}, \, k \ge 4$. \eop \\

We now look at permutations of infinite subsets of integers. Davis et al. 
\cite{Dav77} observed that any permutation of the positive integers 
contains a $3$-term AP as a subsequence. (Let $a_1$ be the first term, and 
let $k$ be the least integer such that $a_k > a_1$. Then $2a_k - a_1$ 
occurs to the right of both $a_1$ and $a_k$.) They also constructed a 
permutation of the positive integers in which no $5$-term AP occurs as a 
subsequence. The corresponding question for $4$-term APs remains open. 
However, if we restrict our attention to arithmetic progressions with odd 
common difference, the problem becomes tractable. \\

\theorem Any permutation of the positive integers must contain a $3$-term 
AP with odd common difference as a subsequence. Furthermore, there exists 
a permutation of the positive integers in which no $4$-term AP with odd 
common difference occurs as a subsequence.

\proof We first show that any 3AP-free permutation 
$\sigma=(t_1,t_2,\ldots,t_{11})$ of $\{1,2,\ldots,11\}$ with $t_1=2$ and 
$t_2=1$ must contain a $3$-term AP with odd common difference as a 
subsequence. Indeed, $4$ must appear in $\sigma$ before $3$, $5$ after 
$4$, $7$ after $4$, $6$ before $7$, $11$ before $6$, and $8$ before $11$. 
Now we have the subsequence $(8,9,10)$ if $9$ occurs before $10$ in 
$\sigma$ and the subsequence $(11,10,9)$ otherwise. This proves our claim. \\

Let $a_1,a_2,\ldots$ be a permutation of the positive integers. Ignoring 
terms less than $a_1$ if necessary, we can assume that $a_1=1$. Let $k$ be 
the least index such that $a_k$ is even, and let 
$a_j=\max(a_1,a_2,\ldots,a_{k-1})$. If $a_j < 2a_k - 1$, then we have 
$(a_1,a_k,2a_k-1)$ as a subsequence. If $a_j \ge 2a_k - 1 > a_k$, let 
$d=a_j - a_k$. Since $a_j$ and $a_j-d$ occur before $a_j+d, 
a_j + 2d, \ldots, a_j + 9d$, and $d$ is odd, it follows from the above 
claim that the permutation contains a $3$-term AP with odd common 
difference. \\

We now exhibit a permutation of the positive integers that contains no 
$4$-term AP with odd common difference as a subsequence. For $i \ge 1$, 
let $\sigma_i$ be a 3AP-free permutation of the following set of $2^i$ 
consecutive even numbers: $$\{(4^i+2)/3, (4^i+8)/3, \ldots, 
(4^{i+1}-4)/3\}$$ Similarly, let $\pi_i$ be a 3AP-free permutation of the 
following set of $2^{i-1}$ consecutive odd numbers: $$\{(4^i+2)/6, 
(4^i+14)/6, \ldots, (4^{i+1}-6)/6\}$$ Observe that the concatenated 
sequence $\sigma_1 \pi_1 \sigma_2 \pi_2 \sigma_3 \pi_3 \cdots$ is a 
permutation of the positive integers. By virtue of our construction, if an 
odd number $x$ occurs in this sequence before an even number $y$, then $2x 
- y < 0$. It follows that no $4$-term AP with odd common difference occurs 
as a subsequence. \eop \\

Following \cite{Sha09}, we will call a subset $S$ of integers $n$-free if 
$S$ can be permuted so that it does not contain any $n$-term AP as a 
subsequence. Given a subset $S$ of the positive integers, let 
$\overline{d}(S)$ and $\underline{d}(S)$ denote, respectively, the upper 
and lower densities of $S$. In other words, $$\overline{d}(S) = \limsup_{n 
\rightarrow \infty} \frac{A(n)}{n} \; \mbox{ and } \; \underline{d}(S) = 
\liminf_{n \rightarrow \infty} \frac{A(n)}{n} \mbox{ where } A(n)=|A \cap 
[1,n]|. $$ Define, for $n \ge 3$, $$\alpha(n)= \sup_S \, 
\{\overline{d}(S): S \mbox{ is $n$-free} \} \; \mbox{ and } \; \beta(n)= 
\sup_S \, \{\underline{d}(S): S \mbox{ is $n$-free} \}. $$ Since the set 
of positive integers is $5$-free, $\alpha(n)=\beta(n)=1$ for $n \ge 5$. 
Bounds for $\alpha(3)$ and $\beta(3)$ were sought in \cite{Dav77}. We 
show the following: \\

\theorem $\alpha(4)=1, \: \alpha(3) \ge 1/2, \: \beta(4) \ge 1/3, 
\: \beta(3) \ge 1/4$. 

\proof Given an integer $a \ge 2$, define $S^{(a)}_i = 
\{a^{2i},a^{2i}+1,\ldots, a^{2i+1}\}$, and let $\sigma^a_i$ be a 
3AP-free permutation of $S^{(a)}_i$. Define $S^{(a)} = \bigcup_{i \ge 0} 
S^{(a)}_i$. We claim that $S^{(a)}$ is $4$-free. Clearly the concatenated 
sequence $\sigma^a_0 \sigma^a_1 \cdots$ does not contain a decreasing 
$3$-term AP. Suppose it contains an increasing $4$-term AP $x_1, x_2, x_3, 
x_4$. Since $x_2, x_3$ and $x_4$ cannot all belong to the same set 
$S^{(a)}_i$, we must have $x_4 \ge 2x_3$ or $x_3 \ge 2x_2$. But then $x_2 
\le 0$ or $x_1 \le 0$, yielding a contradiction. Note that 
$\overline{d}(S^{(a)})=a/(a+1)$ and $\underline{d}(S^{(a)})=1/(a+1)$.  
Since $a$ can be arbitrarily large, it follows that $\alpha(4)=1$. Taking 
$a=2$, we get $\beta(4) \ge 1/3$. \\

Let $p_0=1,q_0=2$, and for $k \ge 1$, define $p_k=2q_{k-1}$ and 
$q_k=3q_{k-1}-1$. Let $\tau_k$ be a 3AP-free permutation of 
$T_k=\{p_k,p_k+1,\ldots,q_k\}$, and let $T=\bigcup_{k \ge 0} T_k$. It is 
easy to verify that $\overline{d}(T)=1/2$ and $\underline{d}(T)=1/4$. We 
claim that the concatenated sequence $\tau_0 \tau_1 \cdots$ contains no 
$3$-term AP as a subsequence. Indeed, if the (increasing) 3AP $x_1, x_2, 
x_3$ occurs as a subsequence, with $x_2$ and $x_3$ belonging to different 
sets $T_k$ and $T_{\ell}$, then $x_3 \ge 2x_2$, so $x_1 \le 0$, yielding a 
contradiction. If $x_2$ and $x_3$ belong to the same set $T_k$, then $x_1 
\in T_{\ell}$ with $\ell < k$. But $x_3 - x_2 < q_k \le x_2 - x_1$, 
contradicting our assumption that $x_1, x_2, x_3$ is a 3AP. Therefore, $T$ 
is $3$-free. Thus $\alpha(3) \ge 1/2$ and $\beta(3) \ge 1/4$. \eop \\

Erd\H{o}s and Graham \cite{Erd80} (see also \cite{Dav77}) asked if it was 
possible to partition the positive integers into two $3$-free sets. 
Clearly, the answer is negative if $\alpha(3) + \beta(3) < 1$. We believe 
this to be the case, and conjecture that the lower bounds in the above 
theorem are optimal, i.e., $\alpha(3)=1/2$ and $\beta(3)=1/4$. However, we 
have not even been able to show that $\beta(3) < 1$.


\begin{thebibliography}{99}

\bibitem{Dav77} J. A. Davis, R. C. Entringer, R. L. Graham and G. J. 
Simmons, {\em {On Permutations Containing No Long Arithmetic 
Progressions}}, Acta Arithmetica 34 (1977), 81-90.

\bibitem{Erd80} P. Erd\H{o}s and R. L. Graham, {\em {Old and New Problems 
and Results in Combinatorial Number Theory}}, L'Enseignment Mathematique, 
Monograph No. 28, Geneva, 1980.

\bibitem{Sha09} A. Sharma, {\em {Enumerating Permutations That Avoid Three 
Term Arithmetic Progressions}}, The Electronic Journal of Combinatorics 16 
(2009), \#R63.

\end{thebibliography}
\end{document}